\documentclass[11pt]{article}
\usepackage{amssymb}
\usepackage[dvips]{graphics}
\usepackage{theorem}

\newtheorem{theorem}{{\sc Theorem}}[section]

\newtheorem{proposition}[theorem]{{\sc Proposition}}
\newtheorem{definition}{{\sc D\'efinition}}[section]

\newtheorem{remarque}{\sc Remarque}
\newtheorem{qquestion}[theorem]{{\sc Question}}
\newtheorem{pprobleme}[theorem]{{\sc Problem}}
\newtheorem{notation}{{\sc Notation}}
\newtheorem{exemple}{\sc Exemple}

\newenvironment{preuve}{{\sc Preuve}}{\nolinebreak
$\Box $}

\newenvironment{demonstration}{{\sc D\'emonstration}} {{\sc CQFD}}

\newenvironment{resume}{\small \begin{center} {\bf R\'esum\'e} \\ } 
{\end{center} \normalsize }

\begin{document}

\title{Un invariant remarquable des polytopes simples}
\author{Fr\'ed\'eric BOSIO}
\maketitle

\begin{resume}
  Nous introduisons un nouvel entier associ\'e \`a tout polytope simple et nous 
\'etudions ses propri\'et\'es.
\end{resume}

\begin{abstract}
  We associate an integer to any simple polytope and we study its properties.
\end{abstract}

\section*{Introduction}

  Dans les ann\'ees 1960, J. Milnor a introduit la notion de rang d'une 
vari\'et\'e (close) $M$ comme le nombre maximum de champs de vecteurs d\'efinis 
dessus, qui commutent, et sont lin\'eairement ind\'ependants en chaque point. 
Le rang d'une vari\'et\'e est tr\`es difficile \`a calculer en g\'en\'eral et 
bien des questions restent pos\'ees \`a son propos. On sait tout de m\^eme que 
le rang de $M$ est strictement positif si et seulement si la caract\'eristique 
d'Euler de $M$ est nulle, ce qui revient \`a demander que $-1$ soit racine du 
polyn\^ome de Poincar\'e $P_M (X)$ de $M$, et une conjecture, attribu\'ee \`a 
Lima, sugg\'erait que pour une vari\'et\'e de rang au moins $2$, $-1$ f\^ut 
aussi racine du polyn\^ome $P'_M $ d\'eriv\'e de $P_M $. Cette conjecture a 
\'et\'e r\'efut\'ee par Bredon~\cite{Br} ; cependant, la valeur en $-1$ de ce 
polyn\^ome, qui n'a gu\`ere fait l'objet d'\'etude jusqu'ici, ne semble pas 
avoir trouv\'e de signification g\'eom\'etrique pr\'ecise.

  Plus r\'ecemment, \`a partir des ann\'ees 1990 s'est d\'evelopp\'ee la 
g\'eom\'etrie torique, au sein de laquelle les vari\'et\'es moment-angle, 
particuli\`erement celles associ\'ees \`a des polytopes simples, occupent une 
place de choix (voir \cite{D-J}, \cite{B-P}). Des liens \'etroits entre 
la combinatoire d'un polytope et la g\'eom\'etrie de la vari\'et\'e 
moment-angle associ\'ee ont notamment \'et\'e \'etablis. G\'en\'eralement, on 
trouve des propri\'et\'es g\'eom\'etriques d'une vari\'et\'e moment-angle 
(rang, nombres de Betti, etc...) \`a partir des propri\'et\'es combinatoires du 
polytope duquel elle provient. Par exemple, le rang d'une telle vari\'et\'e 
vaut au moins la diff\'erence entre le nombre de facettes du polytope et sa 
dimension (on ignore, \`a ma connaissance, s'il y a toujours \'egalit\'e), donc 
au moins $2$ si le polytope consid\'er\'e n'est pas un simplexe.

  Ici, nous allons en quelque sorte dans le sens contraire, en associant un 
invariant entier \`a un polytope \`a partir de la vari\'et\'e moment-angle 
associ\'ee. Nous en calculons la valeur dans les cas les plus communs et 
d\'ecrivons son comportement par rapport \`a quelques unes des op\'erations les 
plus usuelles sur les polytopes simples. Nous terminons par quelques cas 
particuliers qui nous ont paru int\'eressants.

\subsection*{Remerciements}

  L'auteur remercie Laurent Meersseman et le projet COMPLEXE 
(ANR-08-JCJC-0130-01) pour leur soutien financier.

  Les calculs ont \'et\'e effectu\'es \`a l'aide du package LVMB de J. 
Tambour~\cite{Ta} ; qu'il en soit aussi remerci\'e.

\section*{Rappels}

\subsection {Polytopes}

\paragraph{Notations}

  Dans la suite de cet article, on se fixera un polytope simple $P$. En fait, 
l'expression polytope d\'esignera, sauf mention contraire, un polytope 
combinatoire (c'est-\`a-dire que deux polytopes g\'eom\'etriques ayant le 
m\^eme poset de faces seront identifi\'es).

  On notera $d$ la dimension de $P$.

  On notera ${\cal F }$ l'ensemble des facettes (i.e. des ($d-1$)-faces de $P$) 
et $n$ leur nombre. Pour une partie ${\cal I }$ de ${\cal F }$, on note 
$F_{\cal I }$ la r\'eunion des facettes de $F$ qui sont dans ${\cal I }$ (en 
fait, seule la topologie de $F_{\cal I }$ nous int\'eressera, et elle ne 
d\'epend pas de la r\'ealisation de $P$). Il s'agit d'une partie du bord de $P$.

  On notera par ailleurs $K_{\cal I }$ le complexe des faces de $I$, o\`u les 
sommets de $K_{\cal I }$ sont les \'el\'ements de ${\cal I }$ et o\`u un 
ensemble de sommets forme une face de $K_{\cal I }$ si leur intersection  est 
non vide dans $P$ (autrement dit, les faces de $K_{\cal I }$ correspondent aux 
faces de $P$ qui sont intersection de facettes toutes dans ${\cal I }$).

  On montre facilement que $K_{\cal I }$ et $F_{\cal I }$ ont le m\^eme type 
d'homotopie (voir par ex. \cite{B-M}).

\paragraph{Op\'erations classiques}

  Nous d\'ecrivons ici quelques op\'erations classiques que l'on r\'ealise sur 
les polytopes simples.

\begin{description}

\item[Produit]
  Si $Q$ et $R$ sont deux polytopes simples, leur produit cart\'esien est aussi 
un polytope simple. Si $Q$ (resp. $R$) est de dimension $d_Q $ (resp $d_R $) et 
poss\`ede $n_Q $ (resp $n_R $) facettes, alors $Q \times R$ est de dimension 
$d_Q + d_R $ et a $n_Q + n_R $ facettes.

\item[Blending] (voir par ex. \cite{B-P})
  On consid\`ere deux polytopes simples $Q$ et $R$ de m\^eme dimension $d$. On 
choisit un sommet $v_Q $ de $Q$ et un sommet $v_R $ de $R$. On choisit aussi 
une correspondance biunivoque $Q_i \mapsto R_i $ entre les facettes 
$Q_1 ,..., Q_d $ de $Q$ contenant $v_Q $ et les facettes  $R_1 ,..., R_d $ de 
$R$ contenant $v_R $. On appelle alors blending, ou somme connexe, de $Q$ et de 
$R$ le polytope $Q \# R$ obtenu en coupant les sommets $v_Q $ et $v_R $, suivi 
de transformations projectives r\'eelles de $Q_{\mbox{tronqu\'e}}$ et 
$R_{\mbox{tronqu\'e}}$ rendant possible la derni\`ere op\'eration, qui consiste 
\`a "recoller" chaque facette $Q_i $ avec la facette correspondante $R_i $ 
(combinatoirement, c'est toujours r\'ealisable).

\scalebox{.4}{\includegraphics{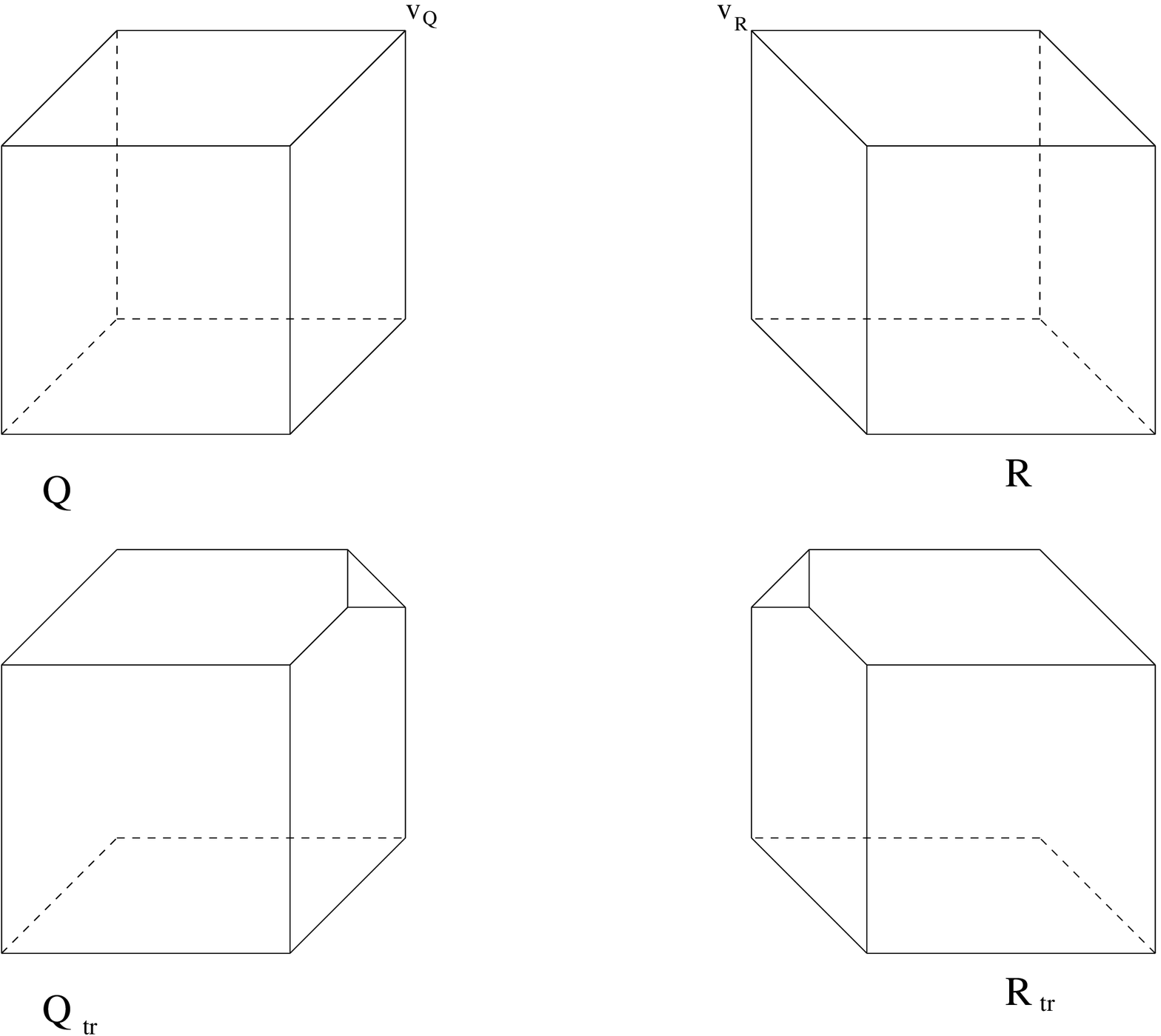}}
\newline
\scalebox{.4}{\includegraphics{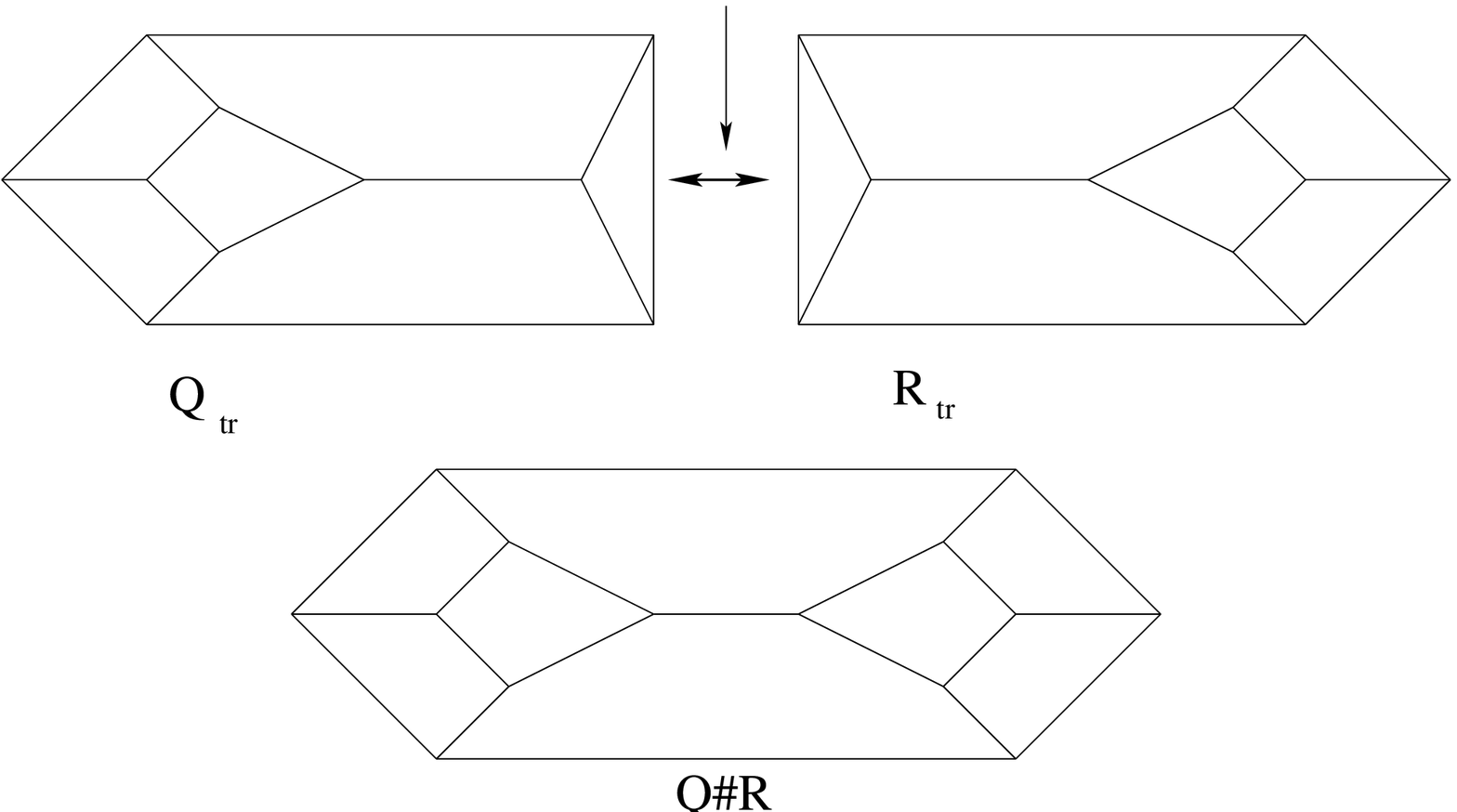}}

  Si $Q$ poss\`ede $n_Q $ facettes et $R$ en poss\`ede $n_R $, alors $Q \# R$ 
en a $n_Q + n_R - d$. Notons tout de m\^eme que le polytope obtenu d\'epend du 
choix des sommets $v_Q $ et $v_R $ ainsi que de la correspondance entre les 
facettes. Toutefois, pour certaines op\'erations, il n'est pas n\'ecessaire de 
les pr\'eciser ; en ce qui concerne l'invariant que nous \'etudions, ce sera 
inutile.

\item[Wedging] (voir par ex. \cite{K-W})
  Si $P$ est un polytope simple et $X$ une facette de $P$, on contruit un 
polytope simple $W_X P$, appel\'e wedge (ou book) de $P$ sur la facette $X$, en 
\'ecrasant, dans le produit de $P$ par un intervalle, la facette $X \times I$ 
sur $X \times \{ 0 \} $.

  On a une projection naturelle de $W_X P $ sur $P$. Les deux facettes de 
$W_X P$ contenant $X \times \{ 0 \} $ seront appel\'ees les facettes 
principales de $W_X P$ ; elles sont combinatoirement \'egales \`a $P$ et leur 
projection r\'ealise cette identification. Toute autre facette se projette 
sur une facette de $P$ autre que $X$ et ces deux facettes se correspondent 
naturellement.

  Si $P$ est de dimension $d$ et a $n$ facettes, alors 
$W_X P$ est de dimension $d+1$ et a $n+1$ facettes. Notons toutefois que le 
polytope obtenu d\'epend de la face $X$ choisie.

\scalebox{.5}{\includegraphics{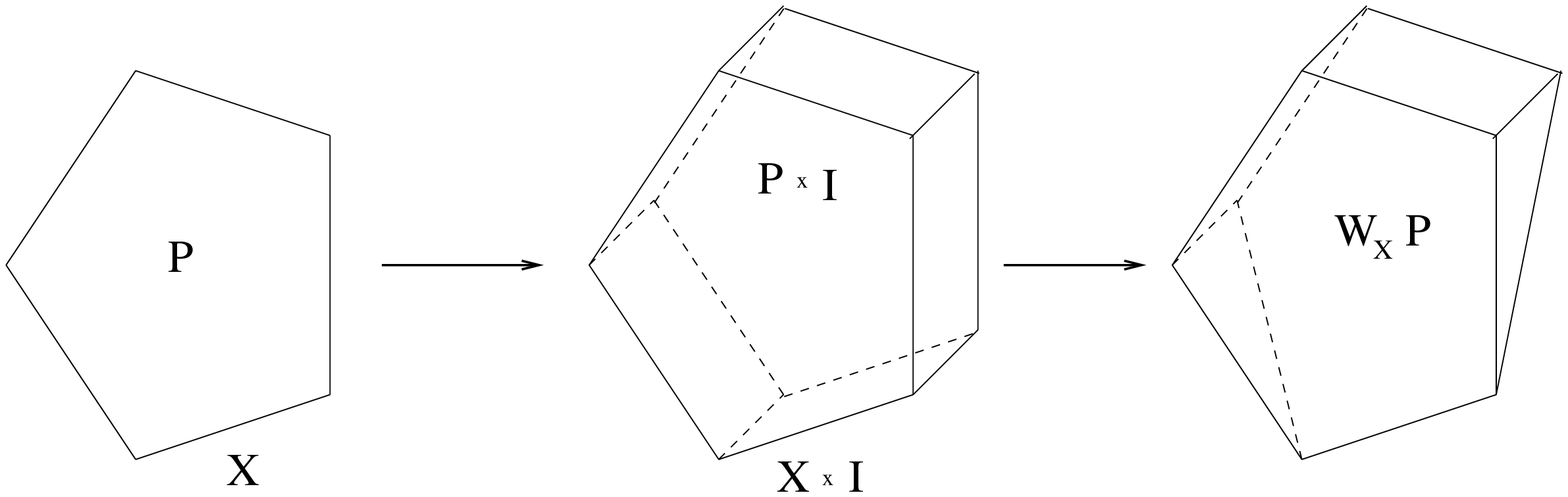}}

\item[Flipping] (voir par ex. \cite{B-M})
  Soit $T$ un polytope simple. On suppose qu'on a deux facettes $Q$ et $R$ de 
$T$ qui sont disjointes et telles qu'il existe un unique sommet $v$ de $T$ qui 
ne soit ni sur $Q$ ni sur $R$. On note alors $a$ (resp. $b$) le nombre 
d'ar\`etes de $T$ issues de $v$ qui arrivent sur $Q$ (resp. sur $R$). On dit 
alors que $R$ est obtenu \`a partir de $Q$ par un $(a,b)$-flip (on a 
$a+b = \dim T$).

  Les deux polytopes $Q$ et $R$ ont m\^eme dimension et, si $a$ et $b$ sont 
$>1$, ils ont m\^eme nombre de facettes.

  La figure suivante pr\'esente un $(2,2)$-flip entre une cube et un livre 
pentagonal~:

\scalebox{.5}{\includegraphics{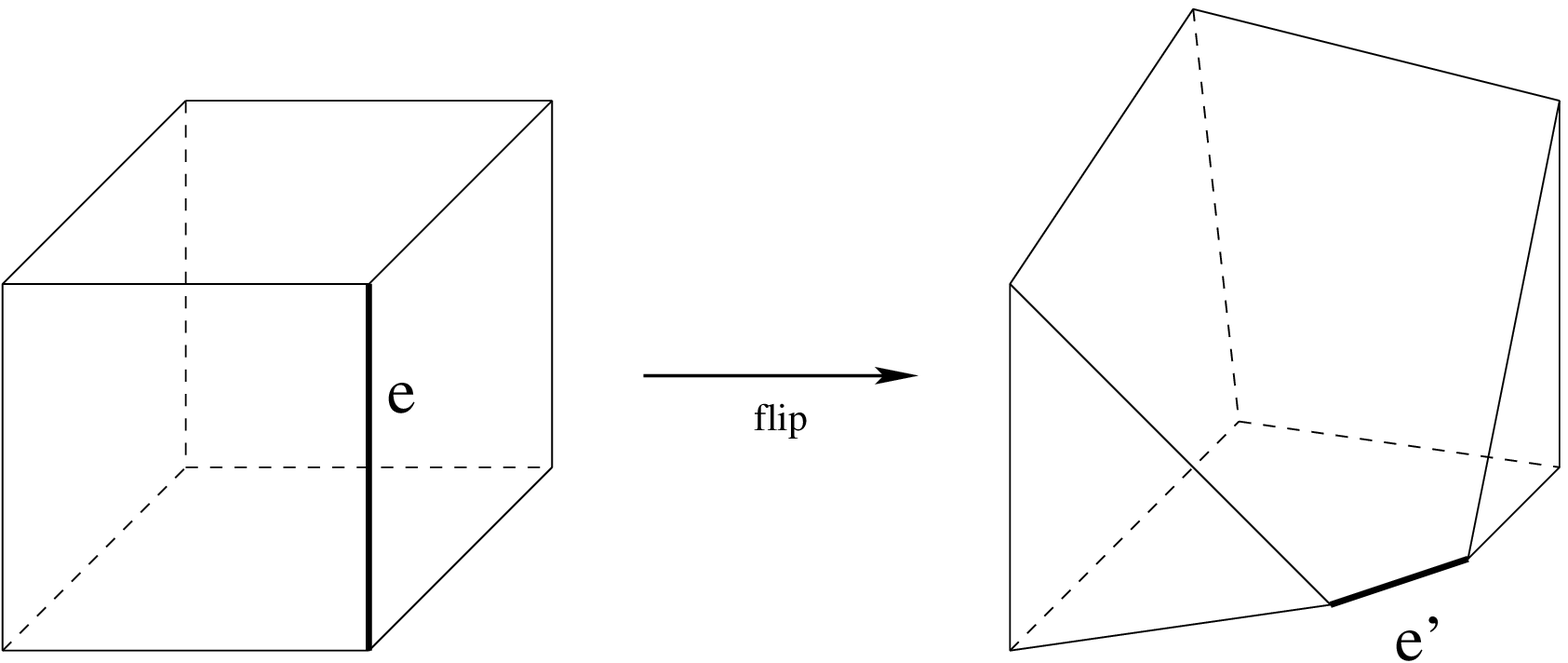}}

\end{description}

\subsection{Vari\'et\'es moment-angle}

\paragraph{Notations}

  Pour un ensemble fini $E$, on note $T_E $ le produit $(S^1 )^E $ o\`u $S^1 $ 
est un cercle, habituellement le cercle unit\'e de $\mathbb C $.

  On note par ailleurs $Z_P $ le complexe moment-angle associ\'e au polytope 
simple $P$ \linebreak
(voir~\cite{B-P}), qui est une vari\'et\'e canoniquement diff\'erentiable. 
Rappelons que $Z_P $ peut \^etre vu comme plong\'e dans $\mathbb C ^{\cal F }$ 
o\`u il est muni de l'action naturelle du tore $T_{\cal F }$ par multiplication 
sur chaque coordonn\'ee.

\paragraph{Homologie}

  En faisant d\'ecrire \`a ${\cal I }$ toutes les parties de ${\cal F }$, on 
peut d\'ecrire l'homologie de $Z_P $ \`a l'aide des $F_{\cal I }$ 
(voir~\cite{B-M}).

  Pour cela, on consid\`ere pour chaque $F_{\cal I }$ son homologie r\'eduite. 
A une classe $[c]$ dans $\tilde{H}_i (F_{\cal I })$, on prend un repr\'esentant 
$c$ qui est un cycle de dimension $i$ dans $F_{\cal I }$. Dans $P$ (vu ici 
comme plong\'e effectivement dans $\mathbb R ^d $), ce cycle est un bord et on 
a donc une cha\^{\i}ne de dimension $i+1$ dont $c$ est le bord. En saturant 
cette cha\^{\i}ne par l'action du tore $T_{\cal I }$, on obtient un cycle de 
dimension $k+i+1$ de $Z_P $, dont la classe ne d\'epend pas des choix 
effectu\'es. En fait, toute l'homologie de $Z_P $ est engendr\'ee par ces 
classes.

N.B. : On consid\`ere plus habituellement le dual de Poincar\'e de cette 
classe, qui est une classe de cohomologie de dimension $n+d-k-i-1$. Il est 
toutefois mieux adapt\'e \`a notre propos de lui associer cette classe 
d'homologie.

\begin{definition}
  La construction ci-dessus nous identifie, pour $R$ anneau quelconque, \newline
$H_* (Z_P ,R)$ avec 
$\bigoplus _{{\cal I } \subset {\cal F }} \tilde{H} (F_{\cal I }, R)$. 

  Si on prend une facette $X$ de $P$, les classes de $H_* (Z_P , R)$ 
correspondant aux \'el\'ements de 
$\bigoplus _{{\cal I } \subset {\cal F } \atop X \in {\cal I}} 
\tilde{H} (F_{\cal I }, R)$, autrement dit celle provenant des ensembles de 
facettes contenant $X$, seront dites $X$-satur\'ees.

  Les classes correspondant au contraire aux \'el\'ements de 
$\bigoplus _{{\cal I } \subset {\cal F } \atop X \notin {\cal I}} 
\tilde{H} (F_{\cal I }, R)$, seront dites $X$-orthogonales.
\end{definition}

\begin{remarque}
\label{dimmoit}
  Par dualit\'e d'Alexander, la dimension de 
$\tilde{H} (F_{\cal I }, {\mathbb R })$ est, pour tout ${\cal I }$, \'egale \`a 
la dimension de $\tilde{H} (F_{{\cal F }\setminus {\cal I }}, \mathbb R )$. Par 
sommation, la dimension de l'espace des classes $X$-satur\'ees est \'egale \`a 
la dimension de l'espace des classes $X$-orthogonales, donc toutes deux valent 
la moiti\'e de la dimension de $H_* (Z_P , {\mathbb R })$.
\end{remarque}

\subsection{Identit\'e calculatoire}

  Nous rappelons ici le r\'esultat suivant, qui nous sera bien utile~:

  Soit $n$ un entier naturel et $P(X)$ un polyn\^ome de degr\'e $<n$. Alors~:

\begin{equation}
\label{ident}
\sum_{k=0}^n (-1)^k P(k) {n \choose k} = 0
\end{equation}

  En effet, en d\'eveloppant $(X+1)^n $, en le d\'erivant $i$ fois 
($0 \leq i < n$) et en prenant la valeur en $-1$, on obtient l'indentit\'e pour 
un polyn\^ome de degr\'e $i$. D'o\`u le r\'esultat par lin\'earit\'e.

  Dans toute la suite pour tout $n \geq 0$ et $k$ ne v\'erifiant pas 
$0 \leq k \leq n$, nous posons ${n \choose k} = 0$. On obtient ainsi aussi 
$\displaystyle{\sum_{k \in \mathbb Z }} (-1)^k P(k) {n \choose k} = 0$.

\section {D\'efinition}

Nous introduisons ici notre nouvel invariant, associ\'e \`a un polytope simple 
$P$ auquel est associ\'ee la vari\'et\'e moment-angle $Z_P $.

  Si $M$ est une vari\'et\'e compacte de dimension $k$, on appelle $b_i (M)$ 
(ou simplement $b_i $ s'il n'y a pas d'ambigu\"{\i}t\'e) son $i$-i\`eme nombre 
de Betti, soit $b_i = \dim H^i (M, \mathbb R )$. On appelle polyn\^ome de 
Poincar\'e de $M$ le polyn\^ome donn\'e par 
$Poinc_M (X) = \sum_{i=0}^{k} b_i X^i $ (on sait que $b_i = 0$ si $i>k$).

  On d\'efinit alors~:

\begin{equation}
I(P) = Poinc'_{Z_P }(-1) = - \sum_{i=0}^{n+d} (-1)^i \cdot i \cdot b_i 
\end{equation}

  Nous allons donner quelques propri\'et\'es de cet invariant (nous verrons par 
rapport \`a quoi il est invariant) en reliant $I(P)$ \`a des propri\'et\'es 
combinatoires de $P$.

\section{Propri\'et\'es de l'invariant}

N.B. : Dans toute la suite, les polytopes consid\'er\'es seront simples sauf 
mention contraire.

  Certaines des valeurs les plus significatives de l'invariant consid\'er\'e 
sont consign\'ees dans le th\'eor\`eme suivant~:

\begin{theorem}
\begin{enumerate}
\item {\bf{Cas d'un simplexe}}

  Notons $\Delta _d $ le $d$-simplexe ($n = d+1$). On a alors 
$I(\Delta _d ) = 2d+1$.

  Ce cas est extr\`emement sp\'ecial (nombre de r\'esultats \`a venir sont faux 
pour lui), par exemple c'est le seul cas o\`u cet invariant est impair ; aussi 
dans la suite l'excluons nous tr\`es fr\'equemment.

\item {\bf{Cas o\`u $n$ et $d$ ont m\^eme parit\'e}}

  On suppose ici que $n$ et $d$ ont m\^eme parit\'e, autrement dit que $n+d$, 
ou de fa\c{c}on \'equivalente $n-d$, est pair. Alors, on a $I(P) = 0$.

\item {\bf{Cas d'un produit}}

  On suppose ici qu'il existe deux polytopes non triviaux (i.e de dimension 
$\geq 1$) $Q$ et $R$ tels que $P = Q \times R$. Alors $I(P) = 0$.

\item {\bf{Cas d'un blending}}

  Supposons ici qu'une facette (au moins) de $P$ soit un $d-1$-simplexe. Alors 
si $n$ et $d$ sont de parit\'e diff\'erentes, et que $P$ n'est pas un simplexe, 
on a $I(P) = 2$.

  Ce cas inclut le cas particulier o\`u $P$ est un polyg\^one \`a un nombre 
impair de c\^ot\'es (triangle except\'e).

  Ce r\'esultat se g\'en\'eralise aux blendings plus g\'en\'eraux. Si on 
suppose qu'il existe deux polytopes $Q$ et $R$ (de dimension $d$) tels que 
$P = Q \# R$, et que $n-d$ est impair (ce qui revient \`a demander que les 
nombres de facettes de $Q$ et de $R$ soient de parit\'es diff\'erentes), on a 
pareillement $I(P) = 2$.

\item {\bf{Cas o\`u $n-d = 3$}}

  On suppose ici que $P$ v\'erifie $n = d+3$. Son dual $P^{*}$ est donc un 
polytope simplicial de dimension $d$ \`a $d+3$ sommets. On lui associe alors 
classiquement un polyg\^one \`a un nombre impair $k$ de c\^ot\'es via son 
diagramme de Gale (voir~\cite{Gr}). On a alors $I(P) = k-3$.

\item {\bf{Cas neighbourly dual}}

  On suppose ici que $P$ est de dimension paire $d = 2d'$ et est neighbourly 
dual, c'est-\`a-dire que si on prend $d'$ facettes quelconques de $P$, leur 
intersection est une $d'$-face de $P$ (non vide suffit en fait car $P$ est 
simple).

  Alors, si $n$ est impair, et que $P$ n'est pas un simplexe, on a $I(P) = d$.

  On peut remarquer que ceci aussi inclut le cas particulier des polyg\^ones 
\`a un nombre impair de c\^ot\'es (triangle except\'e).

\item {\bf{Cas d'un wedge}}

  On suppose ici que $P$, non simplexe, est obtenu par wedge sur un polytope 
$Q$ de dimension $d-1$. On a alors $I(P) = I(Q)$.
\end{enumerate}
\end{theorem}

  Cette derni\`ere propri\'et\'e motive l'appellation d'invariant pour $I$.

\begin{demonstration}
\begin{itemize}

\item Le cas 1), $P = \Delta _d $ est imm\'ediat. Dans ce cas, $Z_{\Delta _d }$ 
est la sph\`ere de dimension $2d+1$ et le r\'esultat en d\'ecoule.

\item Le cas 2) n'est pas difficile non plus. Rappelons que $n+d$ est la 
dimension de $Z_P $. On sait que la caract\'eristique d'Euler de $Z_P $ est 
nulle. De plus, $Z_P $ v\'erifie la dualit\'e de Poincar\'e (car $Z_P $ est une 
vari\'et\'e close et orientable) ; on a donc $b_{n+d-i} = b_i $ pour tout $i$. 
En effectuant le changement d'indice $j = n+d-i $ dans la somme d\'efinissant 
$I(P)$, on obtient $I(P) = -I(P) + 0$. Donc $I(P)$ est bien nul.

\item Le cas 3) est simple. On sait que si $P$ et $Q$ sont deux polytopes, 
alors $Z_{P \times Q} = Z_P \times Z_Q $.

  Si on se donne deux vari\'et\'e $M$ et $N$ et qu'on d\'efinit $I(M)$ et 
$I(N)$ par la m\^eme formule que ci-dessus, la formule de K\"unneth 
(voir par ex. \cite{Ha}) donne 
$I(M \times N) = I(M) \cdot \chi(N) + \chi(M) \cdot I(N)$. Ainsi, si $M$ et $N$ 
ont toutes deux une caract\'eristique d'Euler nulle, alors $I(M \times N) = 0$. 
C'est en particulier le cas pour $Z_P $ et $Z_Q $ d'o\`u $I(P \times Q) = 0$.

\item Le cas 4) est plus compliqu\'e. Il s'inspire de la preuve de la 
proposition~11.2 de \cite{B-M}, mais il faut la g\'en\'eraliser \`a des 
blendings plus g\'en\'eraux.

  Commen\c{c}ons par le cas particulier o\`u $P$ -non simplexe- poss\`ede une 
facette simpliciale. Ceci \'equivaut \`a ce que $P$ soit obtenu en coupant un 
sommet \`a un polytope $Q$, soit encore comme blending $Q \# \Delta _d $.

 On appelle ici $m$ le nombre de facettes de $Q$ et $P$ en a alors $n = m+1$. 
La proposition~11.2 de \cite{B-M} nous dit que pour $3 \leq j \leq m+d-2$, on 
a~:
$$b_i (Z_P ) = b_i (Z_Q ) + b_{i-1} (Z_Q ) + {m-d \choose j-2d+1} + 
{m-d \choose i-2}$$
  D\'esignons $b_i (Z_Q )$ par $b_i $. Le calcul donne alors~:

$$I(P) = (-1)^{m+d} (m+d+1) + 
\sum_{i=3}^{m+d-2} (-1)^{i+1} \cdot i \cdot 
(b_{i-1} + b_i + {m-d \choose i-2d+1} + {m-d \choose i-2}) = $$
$$(-1)^{m+d} (m+d+1) + \sum_{i=3}^{m+d-2} (-1)^{i+1} \cdot i \cdot b_{i-1} + 
\sum_{i=3}^{m+d-2} (-1)^{i+1} \cdot i \cdot b_i + $$
$$\sum_{i=3}^{m+d-2} (-1)^{i+1} \cdot i \cdot {m-d \choose i-2d+1} + 
\sum_{i=3}^{m+d-2} (-1)^{i+1} \cdot i \cdot {m-d \choose i-2}$$
Appelons $C$ la derni\`ere de ces sommes, $B$ l'avant derni\`ere et $A$ la 
somme des deux qui les pr\'ec\`edent.

  D'apr\`es les propri\'et\'es des vari\'et\'es moment-angle, on a 
$b_2 = b_{m+d-2} = 0$. D'o\`u 
$$A = \sum_{i=3}^{m+d-2} (-1)^{i+1} \cdot i \cdot b_{i-1} + 
\sum_{i=3}^{m+d-2} (-1)^{i+1} \cdot i \cdot b_i = $$
$$\sum_{k=3}^{m+d-3} (-1)^k \cdot {k+1} \cdot b_k + 
\sum_{k=3}^{m+d-3} (-1)^{k+1} \cdot k \cdot b_k = 
\sum_{k=3}^{m+d-3} (-1)^k \cdot b_k $$
  La caract\'eristique d'Euler de $Z_Q $ est nulle et, comme 
$b_1 = b_2 = b_{m+d-2} = b_{m+d-1} = 0$ et $b_0 = b_{m+d} = 1$, on trouve 
$0 = A + 1 + (-1)^{m+d}$, soit $A = -1 + (-1)^{m+d+1}$

  D'autre part, quand $i$ varie de $3$ \`a $m+d-2$, alors $i-2d+1$ varie de 
$4-2d \leq 0$ \`a $m-d-1$. D'o\`u 
$$B = \sum_{i=3}^{m+d-2} (-1)^{i+1} \cdot i \cdot {m-d \choose i-2d+1} = 
\sum_{k=0}^{m-d-1} (-1)^k \cdot (k+2d-1) \cdot {m-d \choose k}$$
Or, d'apr\`es l'identit\'e~(\ref{ident}),
$\sum_{k=0}^{m-d} (-1)^k \cdot (k+2d-1) \cdot {m-d \choose k} = 0$, 
donc \newline
$B = 0 - \left( (-1)^{m-d} \cdot 1 \cdot (m+d-1) \right) = 
(-1)^{m+d+1} (m+d-1)$.

  Aussi, $i$ varie de $3$ \`a $m+d-2$, alors $i-2$ varie de 
$1$ \`a $m+d-4 \geq m-d$. D'o\`u 
$$C = \sum_{i=3}^{m+d-2} (-1)^{i+1} \cdot i \cdot {m-d \choose i-2} = 
\sum_{k=1}^{m-d} (-1)^{k+1} \cdot (k+2) \cdot {m-d \choose k}$$
Or, $\sum_{k=0}^{m-d} (-1)^{k+1} \cdot (k+2) \cdot {m-d \choose k} = 0$, d'o\`u 
$C = 0 - (-1) \cdot 2 \cdot 1 = 2$.

  On obtient donc finalement~:
$$I(P) = (-1)^{m+d} (m+d+1) + A + B + C = $$
$$(-1)^{m+d} (m+d+1) + (-1)^{m+d+1} (m+d-1) + (-1 + (-1)^{m+d+1}) + 2 = $$
$$(-1)^{m+d} (m+d+1-m-d+1-1) + (2-1) = 1 + (-1)^{m+d} = 1 - (-1)^{n+d}$$
  Ainsi, si $Z_P $ est de dimension impaire, on a bien $I(P) = 2$.

  Passons maintenant \`a un blending plus g\'en\'eral. Nous allons pour 
commencer fixer un certain nombre de notations.

  On se donne un polytope $P$ de dimension $d$ obtenu comme blending de deux 
polytopes $Q$ et $R$ en des sommets $v_Q $ et $v_R $. Nous verrons que ni les 
sommets choisis ni la mani\`ere dont est effectu\'e le blending ne jouent de 
r\^ole. Nous supposons aussi que ni $Q$ ni $R$ n'est un simplexe (sinon on est 
ramen\'e au cas ci-dessus).

  Notons $n_Q $ (resp. $n_R $) le nombre de facettes de $Q$ (resp. $R$). Notons 
comme d'habitude ${\cal F }$ l'ensemble des facettes de $P$. On note 
${\cal F }_Q $ (resp. ${\cal F }_R $) \`a la fois l'ensemble des facettes de 
$Q$ (resp. $R$) qui ne contiennent pas le sommet en lequel a \'et\'e 
effectu\'e le blending et l'ensemble des facettes de $P$ qui leur 
correspondent. Notons ${\cal F }_v $ \`a la fois l'ensemble des facettes de $P$ 
obtenues par recollement d'une facette de $Q$ et d'une de $R$, ainsi que 
l'ensemble des facettes de $Q$ ou $R$ contenant le sommet o\`u a \'et\'e 
effectu\'e le blending. Ainsi, ${\cal F }$ est partitionn\'e par 
${\cal F }_Q $, ${\cal F }_R $ et ${\cal F }_v $.

  En fait, nous ne distinguerons pas un facette de $Q$ ou de $R$, m\^eme dans 
${\cal F }_v $, de la facette de $P$ qui lui correspond.

  Prenons un ensemble ${\cal I } \subset {\cal F }$ de facettes de $P$.

  On note $F_{\cal I }$ la r\'eunion des facettes de $P$ qui sont dans 
${\cal I }$.

  La d\'ecomposition de ${\cal I }$ dans la partition de ${\cal F }$ ci-dessus 
sera not\'ee ${\cal I } = {\cal I }_Q \cup {\cal I }_R \cup {\cal I }_v $.

ATTENTION : On prendra soin de ne pas confondre $F_{{\cal I }_Q }$, qui est la 
r\'eunion des seules facettes situ\'ees dans ${\cal I }_Q $ avec 
$F_{{\cal I },Q}$, qui est la r\'eunion des facettes de ${\cal I }$ qui 
proviennent de $Q$, autrement dit $F_{{\cal I }_Q \cup {\cal I }_v }$. De 
m\^eme pour $R$.

  Nous allons d\'ecrire l'homologie (r\'eduite) de $F_{\cal I }$.

\begin{itemize}
\item Si ${\cal I }_v $ n'est ni vide ni \'egal \`a tout ${\cal F }_v $, alors 
$F_{\cal I }$ est "la somme connexe de $F_{{\cal I },Q}$ et de 
$F_{{\cal I },R}$ le long de leur fronti\`ere" (les composantes r\'eunies 
\'etant celles contenant ${\cal F }_{{\cal I }_v }$).

 L'homologie r\'eduite de $F_{\cal I }$ est alors la somme directe des 
homologies r\'eduites de $F_{{\cal I },Q}$ et de $F_{{\cal I },R}$.

\item Si ${\cal I }_v $ est vide, alors $F_{\cal I }$ est la r\'eunion 
disjointe de $F_{{\cal I },Q}$ et de $F_{{\cal I },R}$ (qui sont aussi dans ce 
cas ${\cal F }_{{\cal I }_Q }$ et ${\cal F }_{{\cal I }_R }$).

  Si aucun des deux n'est vide, l'homologie r\'eduite de $F_{\cal I }$ peut 
alors \^etre vu comme la somme directe de l'homologie r\'eduite de 
$F_{{\cal I },Q}$, de celle de $F_{{\cal I },R}$ et de $\mathbb Z $ en 
dimension $0$ (non canoniquement plong\'e dedans, un repr\'esentant d'un 
g\'en\'erateur de ce dernier groupe \'etant le cycle form\'e de la diff\'erence 
d'un point de $F_{{\cal I },Q}$ et d'un point de $F_{{\cal I },R}$, la classe 
de ce cycle d\'ependant des composantes des points choisis).

  Si un (seul) des deux est vide $F_{\cal I }$ est soit 
${\cal F }_{{\cal I }_Q }$ et ${\cal F }_{{\cal I }_R }$ et son 
homologie est celle qui correspond. Notons juste que l'homologie r\'eduite de 
l'autre (l'ensemble vide en a en dimension $-1$) n'apparait pas.

\item Si, \`a l'oppos\'e, ${\cal I }_v = {\cal F }_v $, alors $F_{\cal I }$ est 
la somme connexe (en tant que vari\'et\'es de dimension $d-1$) de 
$F_{{\cal I },Q}$ et de $F_{{\cal I },R}$, les endroits de recollement \'etant 
des voisinages de $v_Q $ et $v_R$.

  Si on n'a ni ${\cal I }_Q = {\cal F }_Q $ ni ${\cal I }_R = {\cal F }_R $, 
alors l'homologie r\'eduite de $F_{\cal I }$ est la somme directe de 
l'homologie r\'eduite de $F_{{\cal I },Q}$, de celle de $F_{{\cal I },R}$ et de 
$\mathbb Z $ en dimension $d-2$, la sph\`ere sur laquelle est effectu\'e le 
recollement n'\'etant alors pas un bord.

  Si on a une (seule) des deux \'egalit\'es 
${\cal I }_{Q \mbox{ ou } R} = {\cal F }_{Q \mbox{ ou } R}$, alors 
$F_{\cal I }$ est obtenu en faisant la somme connexe de ${\cal F }_Q $ ou 
${\cal F }_R $ avec une sph\`ere, ce qui ne modifie pas sa topologie et 
l'homologie r\'eduite correspond. Notons juste que l'homologie r\'eduite de 
l'autre (qui est une sph\`ere) n'apparait pas.

\end{itemize}

  En r\'esum\'e, on peut d\'ecomposer l'homologie de $Z_P $ en une homologie 
provenant de $Q$ ou de $R$ et une homologie "sp\'eciale" provenant du 
recollement. Ceci va nous permettre de mener \`a bien notre calcul.

  Prenons par exemple une partie ${\cal J } = {\cal J }_Q \cup {\cal J }_v $ 
des facettes de $Q$. On suppose que ${\cal J }$ n'est ni vide ni tout 
${\cal F }_Q \cup {\cal F }_v  $ et on note $j$ son cardinal. Prenons aussi une 
classe d'homologie r\'eduite $[c]$ de $F_{\cal J }$, disons 
$[c] \in \tilde{H}_l (F_{\cal J }, \mathbb Z )$. Cette classe induit une partie 
de l'homologie de $Z_P $. En fait, pour chaque partie ${\cal K }$ de 
${\cal F }_R $ dont le cardinal est not\'e $k$, elle induit une classe 
$[c]_{\cal K }$ de $H_{j+k+l+1}(Z_P , \mathbb Z )$.

  On a de plus suppos\'e que $R$ n'\'etait pas un simplexe, et donc 
$Card ({\cal F }_R ) \geq 2$, ce qui entra\^{\i}ne que~:
$$\sum_{{\cal K } \subset {\cal F }_R } (-1)^{j+k+l+1} (j+k+l+1) = $$
$$(-1)^{j+l+1} \sum_{k=0}^{Card ({\cal F }_R )} 
(-1)^k {Card ({\cal F }_R ) \choose k} (k + (j+l+1)) = 0$$
  En ce sens, l'homologie induite par $[c]$ n'apporte pas de contribution \`a 
$I(P)$. Ainsi, toute l'homologie provenant de $Q$, et de m\^eme pour celle 
provenant de $R$ peut \^etre occult\'ee dans le calcul de $I(P)$.

  Autrement dit, $I(P)$ r\'esulte uniquement de l'homologie provenant des 
classes sp\'eciales.

  Prenons ${\cal I } \subset {\cal F }$ et notons $i$ le cardinal de 
${\cal I }_Q $, $j$ celui de ${\cal I }_R $.

  Dans le cas o\`u ${\cal I }_v $ est vide, mais ni ${\cal I }_Q $ ni 
${\cal I }_R $, alors "la" classe sp\'eciale associ\'ee dans $Z_P $ est en 
degr\'e $i+j+1$.

  Par dualit\'e, si ${\cal I }_v = {\cal F }_v $ mais que 
${\cal I }_Q \neq {\cal F }_Q $ et ${\cal I }_R \neq {\cal F }_R $, 
alors la classe sp\'eciale associ\'ee dans $Z_P $ est en degr\'e 
$(n_Q + n_R ) - ((Card ({\cal F }_Q ) - i) + (Card ({\cal F }_R ) - j) + 1) = 
i+j+2d-1$.

  Ne reste plus alors \`a prendre en compte que $H_0 (Z_P , \mathbb Z )$ et 
$H_{n_Q + n_R} (Z_P , \mathbb Z )$.

  On obtient donc finalement~:

$$I(P) = (-1)^{n_Q + n_R + 1} (n_Q + n_R ) + $$
$$\sum_{1 \leq i \leq n_Q - d} \sum_{1 \leq j \leq n_R - d}
{n_Q - d \choose i} {n_R - d \choose j} (-1)^{i+j} (i+j+1) + $$
$$\sum_{0 \leq i \leq n_Q - d - 1} \sum_{0 \leq j \leq n_R - d - 1}
{n_Q - d \choose i} {n_R - d \choose j} (-1)^{i+j} (i+j+2d-1)$$

  La premi\`ere des deux sommes peut \^etre d\'ecompos\'ee en trois par 
$$\sum_{1 \leq i \leq n_Q - d} \sum_{1 \leq j \leq n_R - d}
{n_Q - d \choose i} {n_R - d \choose j} (-1)^{i+j} (i+j+1) = $$
$$\left( \sum_{1 \leq i \leq n_Q - d} (-1)^i i {n_Q - d \choose i} \right)
\left( \sum_{1 \leq j \leq n_R - d} (-1)^j {n_R - d \choose j} \right) + $$
$$\left( \sum_{1 \leq i \leq n_Q - d} (-1)^i {n_Q - d \choose i} \right)
\left( \sum_{1 \leq j \leq n_R - d} (-1)^j j {n_R - d \choose j} \right) + $$
$$\left( \sum_{1 \leq i \leq n_Q - d} (-1)^i {n_Q - d \choose i} \right)
\left( \sum_{1 \leq j \leq n_R - d} (-1)^j {n_R - d \choose j} \right) $$

  Or, \`a nouveau d'apr\`es l'identit\'e~(\ref{ident}), on a 
$$\sum_{0 \leq i \leq n_Q - d} (-1)^i {n_Q - d \choose i} = 
\sum_{0 \leq i \leq n_Q - d} (-1)^i i {n_Q - d \choose i} = $$
$$\sum_{0 \leq j \leq n_ - d} (-1)^j {n_R - d \choose j} = 
\sum_{0 \leq j \leq n_ - d} (-1)^j j {n_R - d \choose j} = 0$$

  Ce qui donne donc 
$\sum_{1 \leq i \leq n_Q - d} (-1)^i {n_Q - d \choose i} = 
\sum_{1 \leq j \leq n_R - d} (-1)^i {n_R - d \choose j} = -1$ et \newline
$\sum_{0 \leq i \leq n_Q - d} (-1)^i i {n_Q - d \choose i} =
\sum_{0 \leq j \leq n_ - d} (-1)^j j {n_R - d \choose j} = 0$

  Au total, $\sum_{1 \leq i \leq n_Q - d} \sum_{1 \leq j \leq n_R - d}
{n_Q - d \choose i} {n_R - d \choose j} (-1)^{i+j} (i+j+1) = 0 + 0 + (-1)^2 = 
1$.

  La seconde somme se d\'ecompose aussi en trois par~:
$$\sum_{0 \leq i \leq n_Q - d - 1} \sum_{0 \leq j \leq n_R - d - 1}
{n_Q - d \choose i} {n_R - d \choose j} (-1)^{i+j} (i+j+2d-1) = $$
$$\left( \sum_{0 \leq i \leq n_Q - d - 1} (-1)^i i {n_Q - d \choose i} \right)
\left( \sum_{0 \leq j \leq n_R - d - 1} (-1)^j {n_R - d \choose j} \right) + $$
$$\left( \sum_{0 \leq i \leq n_Q - d - 1} (-1)^i {n_Q - d \choose i} \right)
\left( \sum_{0 \leq j \leq n_R - d - 1} (-1)^j j {n_R - d \choose j} \right) + 
$$
$$(2d-1) \left( \sum_{1 \leq i \leq n_Q - d} (-1)^i {n_Q - d \choose i} \right)
\left( \sum_{1 \leq j \leq n_R - d} (-1)^j {n_R - d \choose j} \right) $$

  D'apr\`es ce qui pr\'ec\`ede,
$\sum_{0 \leq i \leq n_Q - d - 1} (-1)^i i {n_Q - d \choose i} = 
0 - (-1)^{n_Q - d} (n_Q - d)$, \newline
$\sum_{0 \leq j \leq n_R - d - 1} (-1)^j {n_R - d \choose j} = 
-(-1)^{n_R - d}$, 
$\sum_{0 \leq i \leq n_Q - d - 1} (-1)^i {n_Q - d \choose i} = 
-(-1)^{n_R - d}$ et \newline
$\sum_{0 \leq j \leq n_R - d - 1} (-1)^j j {n_R - d \choose j} = 
-(-1)^{n_R - d} (n_R - d)$. D'o\`u~:
$$\sum_{0 \leq i \leq n_Q - d - 1} \sum_{0 \leq j \leq n_R - d - 1}
{n_Q - d \choose i} {n_R - d \choose j} (-1)^{i+j} (i+j+2d-1) = $$
$$(-1)^{n_Q + n_R } (n_Q - d) + (-1)^{n_Q + n_R } (n_R - d) + (2d-1) 
(-1)^{n_Q + n_R } = (-1)^{n_Q + n_R } (n_Q + n_R - 1)$$

  On obtient donc finalement~:
$$I(P) = (-1)^{n_Q + n_R + 1} (n_Q + n_R ) + 1 + 
(-1)^{n_Q + n_R }(n_Q + n_R - 1) = 
1 - (-1)^{n_Q + n_R }$$
  Soit donc bien $2$ si $n_Q + n_R $ est impair.

\item Le cas 5) est en fait un corollaire direct des cas 3), 6) et 7). En 
effet, un polytope simple v\'erifiant $n-d = 3$ est obtenu par des wedges 
successifs au-dessus, soit du cube si $k=3$, le polytope \'etant alors un 
produit de trois simplexes, soit du dual du polytope cyclique $C(k-3,k)$ si 
$k \geq 5$.

  Dans le premier cas, le poit 3) nous assure que $I(P) = 0 = k-3$. Dans 
l'autre cas, le 7) nous dit que $I(P) = I(C^* (k-3,k))$ et le 6) nous dit que 
ceci n'est autre que la dimension de ce dernier polytope, soit $I(P) = k-3$, 
car il est bien connu qu'un polytope cyclique est neighbourly \cite{Gr}.

\item Nous traiterons le cas 6) en calculant explicitement les nombre de Betti 
des vari\'et\'es moment-angle associ\'es au polytopes neighbourly duaux.

\begin{proposition}
  Soit $P$ un polytope neighbourly dual de dimension $d = 2d'$ \`a $n$ facettes 
et soit $j$ tel que $d' < j < n+d'$.

  Alors on a~:

$$b_j (Z_P ) = \frac{1}{(d')!} {n \choose j-d'} 
\prod_{i=0}^{d'-1} \frac{(j - 2d' + i)(n-j+i)}{n-d'+i}$$
\end{proposition}

\begin{preuve}

  Nous allons ici utiliser les complexes $K_{\cal I }$ sur les ensembles de 
facettes de $P$.

  Supposons qie ${\cal I }$ ne soit ni vide ni \'egal \`a ${\cal F }$ tout 
entier. Alors $K_{\cal I}$ a le type d'homotopie d'un bouquet d'un certain 
nombre $b_{\cal I }$ de $(d'-1)$-sph\`eres. Or, si on prend un bouquet de 
$(d'-1)$-sph\`eres, au nombre de $r$, sa caract\'eristique d'Euler vaut 
$1 + (-1)^{d'-1} r$ et la dimension de son $H_{d'-1}$ vaut $r$. D'o\`u 
$b_{\cal I } = (-1)^{d'-1} \cdot (\chi (K_{\cal I }) - 1)$.

  Prenons $d' < j < n+d'$. Alors, une partie ${\cal I }$ de ${\cal F }$ \`a 
$j-d'$ \'el\'ements n'est ni vide ni tout ${\cal F }$ et les classes 
d'homologie qu'elle induit sont de dimension $(j-d') + (d'-1) + 1 = j$. De 
plus, les classes induites par l'ensemble vide ou tout le bord du polytope sont 
en dimension $0$ ou $n + 2d'$, donc pas $j$. On a donc ainsi~:
$$b_j (Z_P ) = 
\sum_{{\cal I } \subset {\cal F } \atop \# {\cal I } = j-d'} b_{\cal I } = 
\sum_{{\cal I } \subset {\cal F } \atop \# {\cal I } = j-d'}
(-1)^{d'-1} \cdot (\chi (K_{\cal I }) - 1) = $$
$$(-1)^{d'} {n \choose j-d'} + 
(-1)^{d'-1} \sum_{{\cal I } \subset {\cal F } \atop \# {\cal I } = j-d'}
\chi (K_{\cal I })$$
  Notons cette derni\`ere somme $b_j $ m\^eme si on ne suppose plus 
$d' < j < n+d'$. Ce n'est alors plus forc\'ement $b_j (Z_P )$.

  Pour calculer $b_j $, on regarde combien de fois chaque face 
de $P$ "apparait" dans les $K_{\cal I }$. Notons tout de m\^eme qu'une $k$-face 
de $P$ est une intersection de $d-k$ facettes et apparait alors comme 
$d-k-1$-simplexe dans les $K_{\cal I }$, et que $d-k-1$ est de parit\'e 
oppos\'ee \`a $k$, donc il faut prendre garde au signe.

  En fait, pour voir appara\^{\i}tre une $k$-face $F$ dans un $K_{\cal I }$, il 
faut et suffit que ${\cal I }$, qui a $j-d'$ \'el\'ements, contienne les $d-k$ 
facettes dont l'intersection est $F$. Il y a alors ${n-d+k \choose j-d'-d+k}$ 
soit ${n-d+k \choose j - 3d' + k}$ complexes $K_{\cal I }$ o\`u cette face 
apparait.

  Cela donne ainsi~:
$$b_j = (-1)^{d'} {n \choose j-d'} + 
(-1)^{d'-1} \sum_{0 \leq k < d} f_k (-1)^{d-k-1}{n-d+k \choose j - 3d' + k} = $$
$$(-1)^{d'} {n \choose j-d'} + 
(-1)^{d'} \sum_{0 \leq k < 2d'} f_k (-1)^k {n - 2d' + k \choose j - 3d' + k}$$
Si on pose $f_{2d'} = 1$, cela donne 
$$b_j = 
(-1)^{d'} \sum_{k=0}^{2d'} (-1)^k f_k {n - 2d' + k \choose j - 3d' + k} = $$
$$(-1)^{d'} {n \choose j-d'}
\left( \sum_{k=0}^{2d'} (-1)^k f_k \prod_{i=0}^{2d' - k - 1} \frac{j-d'-i}{n-i} 
\right) $$
  On constate ainsi que $b_j $ est le produit de ${n \choose j-d'}$ par 
la valeur en $j$ d'un polyn\^ome $P_{n,d}(X)$. De plus, le degr\'e de 
$P_{n,d}(X)$ est major\'e par $2d'$.

  On sait de plus que pour $0 < j \leq 2d'$, on a $b_j (Z_P ) = 0$. Par 
dualit\'e de Poincar\'e, pour $n \leq j < n+2d'$, on a aussi $b_j (Z_P ) = 0$. 
Or, pour $d' < j \leq 2d'$, on a 
$0 = b_j (Z_P ) = b_j = (-1)^{d'} {n \choose j-d'} P_{n,d}(j)$. Comme alors 
${n \choose j-d'}$ n'est pas nul, c'est que $P_{n,d}(j) = 0$.

  Le m\^eme raisonnement montre aussi que $P_{n,d}(j) = 0$ si $n \leq j < n+d'$.

  Ainsi, on a trouv\'e $2d'$ racines au polyn\^ome $P_{n,d}(X)$. Vu son 
degr\'e, c'est un multiple de $\prod_{i = 0}^{d'-1} (X - 2d' + i)(X - (n+i))$.

  Pour trouver lequel, calculons $P_{n,d}(d')$. On a \newline
$P_{n,d}(d') = (-1)^{d'} \sum_{k=0}^{2d'} 
(-1)^k f_k \prod_{i=0}^{2d' - k - 1} \frac{-i}{n-i}$. Si $k < 2d'$, le produit 
pr\'ec\'edent contient le facteur $-i$ pour $i=0$, donc est nul. Il ne reste 
donc que $P_{n,d}(d') = (-1)^{d'} (-1)^{2d'} f_{2d'} \cdot 1 = (-1)^{d'} f_d $ 
et comme $f_d = 1$, on a $P_{n,d}(d') = (-1)^{d'}$.

  Or, posons $Q_{n,d}(X) = \frac{1}{(d')!}
\prod_{i=0}^{d'-1} \frac{(X - 2d' + i)(n-X+i)}{n-d'+i}$. Il est aussi de 
degr\'e $2d'$, s'annule pour les entiers $i$ tels que $d' < i \leq 2d'$ ou 
$n \leq i < n+d'$, et~:
$$Q_{n,d}(d') = 
\frac{1}{(d')!} \prod_{i=0}^{d'-1} \frac{(i-d')(n-d'+i)}{n-d'+i} = 
(-1)^{d'} \frac{1}{(d')!} \prod_{i=0}^{d'-1} (d'-i) = (-1)^{d'} = P_{n,d}(d')$$

  Les polyn\^omes $P_{n,d}(X)$ et $Q_{n,d}(X)$ sont donc \'egaux, d'o\`u la 
formule.
\end{preuve}

  Cette proposition permet alors de calculer $I(P)$. Posons, pour tout entier 
relatif $j$, $\tilde{b}_j = \frac{1}{(d')!} {n \choose j-d'} 
\prod_{i=0}^{d'-1} \frac{(j-d+i)(n-j+i)}{n-d'+i}$. Alors, la formule 
$b_{j} (Z_P ) = \tilde{b}_j $ est en fait valable pour tous les $j$, except\'e 
les suivants~:

\noindent
i) $j=0$, $j = n+d$ : $\tilde{b}_j  = 0$ tandis que $b_j (Z_P ) = 1$ ; \newline
ii) $j = d'$, $j=n+d'$ : $\tilde{b}_j  = (-1)^{d'}$ tandis que $b_j (Z_P ) = 0$.

  On a de plus $\sum_{j \in \mathbb Z } (-1)^j j \tilde{b}_j = 0$. En effet, en 
posant $j' = j-d'$, on constate que $\tilde{b}_j $ est de la forme 
$Q(j') {n \choose j'}$ o\`u $Q$ est un polyn\^ome de degr\'e $d$. D'o\`u la 
somme voulue \'egale au signe pr\`es \`a 
$\sum_{j' \in \mathbb Z } (-1)^{j'} [(X+d')Q](j') {n \choose j'}$. Comme 
$n > d+1$ car le cas du simplexe est exclus, on a bien 
$\sum_{j \in \mathbb Z } (-1)^j \tilde{b}_j = 0$.

  Finalement, pour calculer $I(P)$, il suffit de calculer sa diff\'erence avec 
cette somme. D'apr\`es ce qui pr\'ec\`ede, et comme $d$ est pair, on a 
$$I(P) = $$
$$0 (1-0) - (-1)^{d'} d' (0-(-1)^{d'}) - 
(-1)^{n+d'} (n+d')(0-(-1)^{d'}) - (-1)^{n+d} (n+d) (1-0) = $$
$$0 + d' + (-1)^n (n+d') - (-1)^{n} (n+d) = d' + (-1)^{n+1} d'$$
  Cela donne bien $I(P) = d$ si $n$ est impair.

\item Pour le 7), nous remarquons que cela \'equivaut en un sens \`a ce qu'une 
facette d'un polytope (autre que le simplexe) induise "autant d'homologie en 
dimension paire qu'impaire". Nous relierons cela au quotient de $Z_P $ par 
l'action naturelle du cercle $S^1 $ correspondant \`a note facette.

  Nous n'allons consid\'erer ici que la partie sans torsion de l'homologie des 
espaces consid\'er\'es, la torsion n'influen\c{c}ant pas le calcul de 
$I(W_X P)$.

  Posons ${\cal F }_W $ l'ensemble des facettes de $W_X P$, ${\cal F }' $ 
l'ensemble ${\cal F } \setminus \{ X \} $ et ${\cal F }'_W $ l'ensemble des 
facettes non principales de $W_X P$. On a d\'ej\`a rappel\'e que les 
\'el\'ements de ${\cal F }' $ et de ${\cal F }'_W $ se correspondent 
naturellement.

  Consid\'erons une partie ${\cal J }$ de ${\cal F }'$ dont on note $k$ le 
cardinal et la partie ${\cal J }_-$ de ${\cal F }'_W $ qui lui correspond. Il 
y a quatre parties de ${\cal F }_W $ dont l'intersection avec ${\cal F }'_W $ 
est ${\cal J }_- $, \`a savoir ${\cal J }_- $ elle-m\^eme, ${\cal J }_1 $ et 
${\cal J }_2 $ qui contiennent chacune une des deux facettes principales 
$F_1 $ et $F_2 $ et ${\cal J }_+ $ qui les contient toutes deux.

  Un ensemble $F_{{\cal J }_i }, i = 1,2$ se r\'etracte par d\'eformation sur 
$F_i $, donc est contractile, donc sans homologie r\'eduite.

  $F_{{\cal J }_- }$ se r\'etracte par d\'eformation sur sa projection sur $P$, 
donc a le m\^eme type d'homoto\-pie, et a fortiori la m\^eme homologie 
r\'eduite, que $F_{\cal J }$.

  Une classe de $\tilde{H}_i (F_{\cal J })$ et la classe de 
$\tilde{H}_i (F_{{\cal J }_- })$ qui lui correspond induisent chacune une 
classe d'homologie de degr\'e $i+k+1$, respectivement de $Z_{W_X P}$ et de 
$Z_P $.

  Consid\'erons maintenant $F_{{\cal J }_+ }$. C'est, dans le bord de $W_X P$, 
le compl\'ementaire de $F_{({\cal F }' \setminus {\cal J })_- }$ (plus 
exactement de sons int\'erieur mais ils ont m\^eme type d'homotopie). Par 
dualit\'e d'Alexander, on a pour tout entier $i$, 
$b_i (F_{{\cal J }_+ }) = b_{n+d-i} (F_{({\cal F }' \setminus {\cal J })_- })$. 
Or, d'apr\`es ce qui pr\'ec\`ede, c'est aussi \'egal \`a 
$b_{n+d-i-1} (F_{{\cal F }' \setminus {\cal J }})$ qui, \`a nouveau par 
dualit\'e d'Alexander, vaut aussi $b_{i-1} (F_{{\cal J } \cup \{ X \} })$.

  Une classe de $\tilde{H}_{i-1} (F_{{\cal J } \cup \{ X \} })$ induit une 
classe d'homologie de degr\'e $(i-1) + (k+1) + 1 = i+k+1$ de $Z_P $, tandis 
qu'une classe de $\tilde{H}_i (F_{{\cal J }_+ })$ induit une classe d'homologie 
de degr\'e $i + (k+2) + 1 = i+k+3$ de $Z_X P$.

  Ainsi, dans un sens, une classe d'homologie de degr\'e $m$ de $Z_P $ induit 
une classe d'homolo\-gie de degr\'e $m+2$ de $Z_{W_X P }$ si elle est 
$X$-satur\'ee et de degr\'e $m$ si elle est $X$-orthogonale.

  Pour tout $j$, appelons $b^s _j  $ la dimension de l'espace des classes 
$X$-satur\'ees de degr\'e $j$ et $b^{ns} _j $ la dimension de l'espace des 
classes $X$-orthogonales.

  On a donc $b_j (Z_P ) = b^s _j + b^{ns} _j $ et 
$b_j (Z_{W_X P}) = b^s _{j-2} + b^{ns} _j $.

  Or, on peut relier la notion de $X$-saturation au quotient $Z_P / T_X $ de 
$Z_P $ par l'action du cercle $T_X $. On a en effet la proposition suivante~:

\begin{proposition}
  On consid\`ere la surjection canonique $\pi $ de $Z_P $ sur $Z_P / T_X $, 
l'appli\-ca\-tion $\pi ^* $ de $H^* (Z_P / T_X ,\mathbb R )$ dans 
$H^* (Z_P ,\mathbb R )$ et on appelle $A$ l'image de $\pi ^* $.

  Alors $\pi ^* $ est injective et on a l'\'equivalence entre~: \newline
i) $c$ est dans $A$ ; \newline
ii) le dual de Poincar\'e de $c$ est une classe $X$-satur\'ee ; \newline
iii) $c$ s'annule sur toutes les classes $X$-satur\'ees.
\end{proposition}

\begin{preuve}
  L'injectivit\'e de l'application r\'esulte directement du fait que 
$Z_P / T_X $ se plonge dans $Z_P $ en tant que retract. Ceci peut, par exemple, 
se voir en remarquant que l'action naturelle de $S^1 $ sur $\mathbb C ^n $ par 
rotation sur la premi\`ere coordonn\'ee a un espace d'orbites qui s'identifie 
\`a $\mathbb R _+ \times \mathbb C ^{n-1}$, donc est un retract de 
$\mathbb C ^n $, et l'action sur $Z_P $ en est une sous-action.

  Prenons une classe $c$ dans $A$, $c = \pi ^* (c')$ o\`u 
$c' \in H^k (Z_P / T_X ,\mathbb R )$ ainsi qu'un cycle satur\'e 
$\alpha \in H_k (Z_P ,\mathbb R )$. On a alors 
$c(\alpha ) = c' (\pi _* \, \alpha )$ et $\pi _* \, \alpha $ est "de dimension 
$(k-1)$", d'o\`u $c' (\pi _* \alpha ) = 0$.

  Soit maintenant une classe de cohomologie dont le dual est $X$-satur\'e. Si 
on l'applique \`a une autre classe $X$-satur\'ee, on obtient l'intersection 
(num\'erique) de deux classes $X$-satur\'ees. Or, l'intersection homologique de 
deux classes $X$-satur\'ees \'etant elle-m\^eme $X$-satur\'ee, leur 
intersection num\'erique est nulle.

  Pour voir maintenant qu'un cocycle qui s'annule sur toute classe 
$X$-satur\'ee est dans $A$ et est le dual d'un cycle satur\'e, on peut se baser 
sur la dimension de ces espaces.

  La dimension de $A$ est la moiti\'e de la dimension de 
$H^* (Z_P ,\mathbb R )$. De m\^eme, d'apr\`es la remarque~\ref{dimmoit}, la 
dimension de l'espace des classes $X$-satur\'ees est aussi la moiti\'e de la 
dimension de $H^* (Z_P ,\mathbb R )$. C'est \'egalement la dimension de 
l'espace des classes qui s'annulent sur toutes les classes $X$-satur\'ees.

  On en d\'eduit que ces trois espaces sont les m\^emes.
\end{preuve}

  Cette proposition nous permet d'affirmer que la dimension de l'espace des 
classes $X$-satur\'ees de degr\'e $j$ de $Z_P $ est \'egale \`a la dimension 
des classes de degr\'e $n+d-j$ de $A$, donc aussi \`a la dimension 
des classes de degr\'e $n+d-j$ de $H^* (Z_P / T_X , \mathbb R )$, soit 
\newline
$b^s _j = b_{n+d-j}(Z_P / T_X )$.

  On obtient alors~:
$$b_j (Z_{W_X P}) = b^{ns} _j + b^s _{j-2} = 
b_j (Z_P ) - b^s _j + b_{n+d-j+2}(Z_P / T_X ) = $$
$$b_j (Z_P ) - b_{n+d-j}(Z_P / T_X ) + b_{n+d-j+2}(Z_P / T_X )$$
Et donc~:
$$I(W_X P) = - \sum_{j=0}^{n+d+2} (-1)^j j \cdot b_j  (Z_{W_X P}) = $$
$$- \sum_{j=0}^{n+d+2} (-1)^j j \cdot b_j  (Z_P) + 
\sum_{j=0}^{n+d+2} (-1)^j j \cdot b_{n+d-j}(Z_P / T_X ) - $$
$$\sum_{j=0}^{n+d+2} (-1)^j j \cdot b_{n+d-j+2}(Z_P / T_X ) = $$
$$I(P) + \sum_{j'=-2}^{n+d} (-1)^{n+d-j'} (n+d-j') \cdot b_{j'}(Z_P / T_X ) - $$
$$\sum_{j'=0}^{n+d+2} (-1)^{n+d-j'+2} (n+d-j'+2) \cdot b_{j'}(Z_P / T_X ) = $$
$$I(P) + \sum_{j'=0}^{n+d} (-1)^{n+d-j'} (-2) \cdot b_{j'}(Z_P / T_X ) = 
I(P) + 2 (-1)^{n+d+1} \sum_{i=0}^{n+d} (-1)^i b_i (Z_P / T_X ) = $$
$$I(P) + 2 (-1)^{n+d+1} \chi (Z_P / T_X )$$

  Or, si $P$ n'est pas un simplexe, le cercle diagonal du tore  
$T_{{\cal F } \setminus \{ X \} }$ agit librement sur $Z_P / T_X $, ce qui 
implique $\chi (Z_P / T_X ) = 0$.

  On obtient donc bien $I(W_X P) = I(P)$.
\end{itemize}
\end{demonstration}

\section{Exemples divers}

  Nous donnons ici diff\'erents exemples qui illustrent aussi le comportement 
de notre invariant.

\begin{exemple}
  Une question naturelle au vu des propri\'et\'es pr\'ec\'edentes est de savoir 
si on a $I(P) \geq 0$ pour tous les polytopes simples. Ce n'est pas le cas, et 
nous pr\'esentons deux contre-exemples. Le premier avec $d=3$, le second avec 
$n-d = 5$, ce qui est dans chaque cas la valeur minimale.

  Soit $P$ l'octa\`edre tronqu\'e, ou permuta\`edre de dimension $3$ (il 
poss\`ede $14$ facettes). Le calcul de l'homologie de $Z_P $ nous donne 
$I(P) = -2$.

  Soit $Q$ le polytope obtenu de la fa\c{c}on suivante~: On place un hypercube 
de dimension $4$ en position de "diamant", i.e. avec un sommet en haut et un 
sommet en bas. On le coupe "horizontalement" par un hyperplan $H$ situ\'e juste 
au-dessous de son hyperplan milieu (i.e. l'hyperplan m\'ediateur des deux 
sommets susmentionn\'es). Le polytope $Q$ est alors la partie sup\'erieure de 
l'hypercube (celle situ\'ee au-dessus de $H$). On constate facilement que $H$ 
sectionne toutes les facettes de l'hypercube, ce qui fait que $Q$ poss\`ede $9$ 
facettes.

  Le calcul de l'homologie de $Z_Q $ nous fournit $I(Q) = -4$.
\end{exemple}

\begin{exemple}
  Une remarque sur $I(P)$, qui est donn\'e comme forme lin\'eaire sur 
l'ensemble des nombres de Betti d'une vari\'et\'e, est que si on prend une 
suite (forc\'ement finie) de polytopes $P_0 ,..., P_k $ dont les nombres de 
Betti de chaque degr\'e sont en progression arithm\'etique, alors leurs 
invariants $I_0 ,..., I_k $ sont aussi en progression arithm\'etique.

  Un exemple peut \^etre fourni par des polytopes obtenus \`a partir 
d'hypercubes en flippant des ar\^etes. Plus pr\'ecis\'ement, prenons 
$d \geq 4 $ et posons $P_0 $ l'hypercube de dimension $d$, qui a $2d$ facettes. 
Il est possible de construire des polytopes $P_1 ,..., P_d $ tels que chaque 
$P_i , i \geq 1$ est obtenu \`a partir de $P_{i-1}$ en flippant une ar\^ete. 
On peut voir que les nombres de Betti des $Z_{P_i }$ ne d\'ependent pas des 
ar\^etes flipp\'ees, et sont en progression arithm\'etique \cite{Bo}. 
(Attention, s'il est n\'ecessaires que deux ar\`etes flipp\'ees quelconques 
soient disjointes et non parall\`eles, on n'obtient pas forc\'ement un polytope 
combinatoire en flippant n'importe lesquelles v\'erifiant cela). Leurs 
invariants sont donc aussi en progression arithm\'etique. En fait, on peut voir 
que si $d$ est impair, alors $I_i = 2i$ pour tout $0 \leq i \leq d$ (si $d$ est 
pair $I_i $ est forc\'ement nul).

  D'autres exemples peuvent \^etre obtenus en flippant des ar\`etes \`a partir 
du polytope $Q$ ci-dessus, ou en flippant d'autres simplexes \`a partir de 
produits, etc...
\end{exemple}

\begin{exemple}
  On peut d'ailleurs se demander s'il est possible d'\'etablir un rapport entre 
l'invariant d'un polytope $P$ et celui d'un polytope obtenu \`a partir de $P$ 
par une transformation \'el\'ementaire (i.e. un unique flip).

  Dans le cas d'un $(1,d)$-flip, on rajoute ou enl\`eve une facette et il est 
donc forc\'e que l'un des deux polytopes ait un invariant nul, quel que soit 
la valeur de l'invariant de l'autre.

  Dans d'autres cas de flips, les invariants ont-ils forc\'ement des valeurs 
proches ?

  Prenons l'exemple du dod\'eca\`edre $D$. Le calcul nous fournit $I(D) = 12$. 
Appelons $D'$ le polytope obtenu \`a partir de $D$ en flippant une ar\`ete 
($D'$ ne d\'epend pas de l'ar\`ete flipp\'ee). Le calcul nous fournit 
$I(D') = 6$. Quelles valeurs obtient-on apr\`es le flip d'une ar\`ete de $D'$ 
(remarquons que toutes peuvent le subir) ?

  On sait d\'eja qu'on obtiendra $12$ en flippant celle qui vient 
d'appara\^{\i}tre (on revient au dod\'eca\`edre), et qu'on obtiendra $2$ si on 
flippe une ar\`ete d'un quadrilat\`ere (car celui-ci sera transform\'e en 
triangle \`a l'issue du flip). Le diagramme suivant r\'ecapitule les valeurs 
obtenues apr\`es flips des diff\'erentes ar\`etes de $D'$~:

\scalebox{.5}{\includegraphics{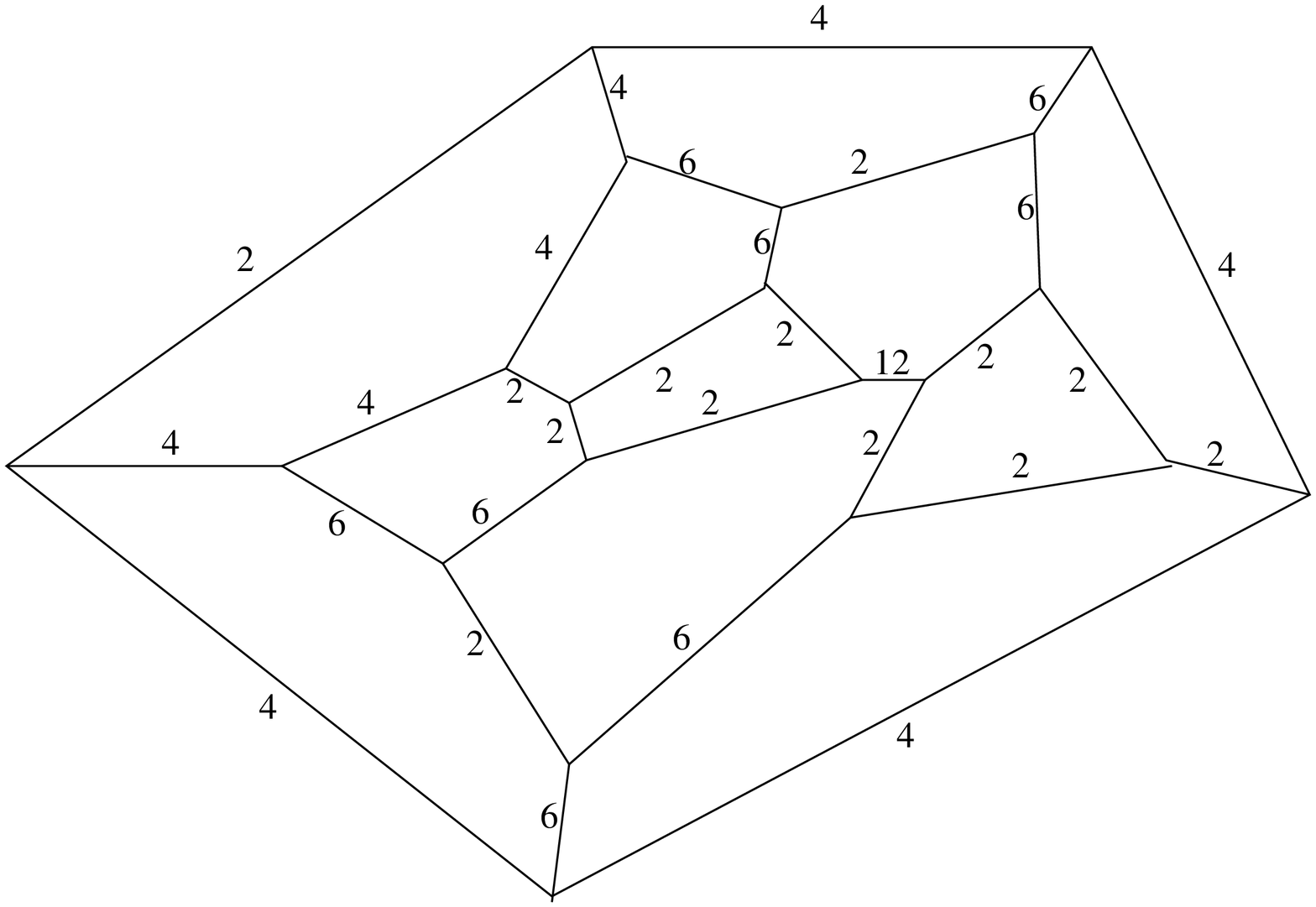}}

  On constate qu'en dehors du cas tr\`es particulier o\`u on revient au 
dod\'eca\`edre, l'invariant ne d\'epasse pas $6$. Ceci invite \`a penser que 
"en r\`egle g\'en\'erale, un poly\`edre (limitons-nous \`a la dimension $3$) 
aura un invariant positif mais assez petit", et nous incite \`a t\'emoigner un 
int\'er\^et particulier aux exceptions.
\end{exemple}

\section{Questions}

  Nous terminons notre propos par quelques questions concernant $I(P)$.

  La premi\`ere, et sans doute la plus importante, est la d\'ependance ou non 
de $I(P)$ par rapport au corps de base. Si $I(P)$ n'en d\'epend pas, 
l'appellation d'invariant pour $I$ est encore plus justifi\'ee, et cela 
renforce encore l'intimit\'e du lien entre $I(P)$ et la combinatoire de $P$.

  En fait, demander que $I(P)$ soit ind\'ependant du corps de base revient \`a 
demander que "pour tout nombre $p$ premier, l'homologie de $Z_P $ (\`a 
coefficients dans $\mathbb Z $) ait autant de $p$-torsion en dimension paire 
qu'en dimension impaire". Or on sait actuellement tr\`es peu de choses 
concernant la torsion dans l'homologie des vari\'et\'es moment-angle.

  Une autre question naturelle est de trouver quelles peuvent \^etre les 
valeurs que prend cet invariant, notamment en fonction de $n$ et $d$. D'abord, 
est-il born\'e ou non lorsque $d$ est fix\'e ? C'est le cas pour $d = 2$ mais 
on peut avoir tendance \`a penser le contraire pour $d \geq 3$ ; toutefois, on 
ignore la r\'eponse exacte. Autre question, a-t-on dans tous les cas 
$|I(P)| \leq n+d$ ? Remarquons que, pour de simples raisons de parit\'e, 
l'\'egalit\'e n'a lieu que dans le cas des simplexes. On peut aussi chercher 
\`a d\'ecrire les polytopes qui r\'ealisent les extrema de $I(P)$ en foction de 
$n$ et $d$.

  On peut aussi chercher des bornes sur $I(P)$ en fonction du rang de Betti 
total $r(Z_P )$ de la vari\'et\'e (i.e. la dimension de 
$H^* (Z_P , \mathbb R )$, soit encore la somme de ses nombres de Betti). Par 
exemple, pour $n-d = 3$, on a $I(P) = \frac{1}{2} r(Z_P ) - 4$. Il semble 
tr\`es raisonnable de conjecturer qu'on a $|I(P)| \leq \frac{1}{2} r(Z_P ) - 4$ 
pour tout polytope simple autre qu'un simplexe.

\section*{Conclusion}

  Toutes ces propri\'et\'es attestent bien de l'int\'er\^et de l'invariant 
$I(P)$. Il semble remarquable qu'en d\'epit de sa d\'efiniton plut\^ot 
g\'eom\'etrique, il soit intimement li\'e \`a la combinatoire du polytope $P$, 
et il est frappant qu'il ait un comportement aussi simple vis-\`a-vis des 
op\'erations usuelles sur les polytopes simples. Toutefois, on aurait fort 
envie de lui trouver une signification plus pr\'ecise de nature combinatoire, 
g\'eom\'etrique ou autre, et, dans le m\^eme ordre d'id\'ees, une m\'ethode qui 
puisse permettre de calculer sa valeur pour des polytopes assez g\'en\'eraux 
auxquels la vari\'et\'e moment-angle associ\'ee a des nombres de Betti dont le 
calcul semble hors de port\'ee.

  En tout cas, il y a fort \`a parier que ce nouvel invariant n'a pas encore 
livr\'e tous ses secrets.

\vskip 10mm

{\footnotesize {Bosio Fr\'ed\'eric \\
Universit\'e de Poitiers \\
UFR Sciences SP2MI \\
D\'epartement de Math\'ematiques \\
UMR CNRS 6086 \\
Teleport 2 \\
Boulevard Marie et Pierre Curie \\
BP 30179 \\
86962 Futuroscope Chasseneuil CEDEX 

e-mail~: bosio@math.univ-poitiers.fr}}

\end{document}